\documentclass[12pt]{article}

\usepackage{amsmath, amssymb, graphicx, amsthm, relsize}
\usepackage{authblk}

\newtheorem{theorem}{Theorem}[section]
\newtheorem{lemma}[theorem]{Lemma}

\newcommand{\Irrc}[1]{\rm Irr(\textit G)}

\newcommand{\acdc}[1]{\rm acd(\textit G)}

\newcommand*\sfrac[2]{{}^{#1}\!/_{#2}}

\begin{document}
\title{\textbf{The Average of Some Irreducible Character Degrees \footnote{This work was done while the first author was a PhD student at Kent State University and this paper appears as a part of his PhD dissertation under the supervision of the second author.}}}
\author{\textit{Ramadan Elsharif}} 
\affil{\textit{Kent State University}  \\ \textit{address: 1300 Lefton Esplanade, Department of Mathematical Sciences, Kent State University, Kent,  OH 44242} \\ \textit{email: relshari@kent.edu}}
\author{\textit{ Mark L. Lewis}}
\affil{\textit{Kent State University}  \\ \textit{address: 1300 Lefton Esplanade, Department of Mathematical Sciences, Kent State University, Kent,  OH 44242} \\ \textit{email: lewis@math.kent.edu}}

\date{}

\maketitle
\vspace{-0.7cm}
\begin{abstract} 
	We are interested in determining the bound of the average of the degrees of the irreducible characters whose degrees are not divisible by some prime $p$ that guarantees a finite group $G$ of odd order is $p$-nilpotent. We find a bound that depends on the prime $p$.  If we further restrict our average by fixing a subfield $k$ of the complex numbers and then compute the average of the degrees of the irreducible characters whose degrees are not divisible by $p$ and have values in $k$, then we will see that we obtain a bound that depends on both $p$ and $k$. Moreover, we find examples that make those bounds best possible. 
	
	
	
\end{abstract}

\newpage

\section{Introduction}
In this paper, all groups are finite, and for a group $G$, ${\rm Irr}(G)$ is the set of (complex) irreducible characters of $G$. The average of the (complex) irreducible character degrees of $G$ is  
$$\displaystyle {\rm acd}(G)= \frac{\sum_{\chi \in {\rm Irr(G)}}^{} \chi(1)}{|{\rm Irr}(G)|}.$$

The influence of character degrees on the structure of finite groups has been a subject of interest for many scholars. For instance, the average of the character degrees can determine the solvability of a group. Indeed, there is interest in determining the best bound for  the average of the irreducible character degrees that will guarantee solvability, nilpotency, or supersolvability of the group.

The story of the average of irreducible character degrees starts when Magaard and Tong-Viet prove in Theorem 1.4 of \cite{MR2764923} that if ${\rm acd}(G)\le2$, then  $G$ is solvable. They conjecture that this bound should be $3$. In the same direction, Isaacs, Loukaki, and Moret\'o in Theorem~A of \cite{MR3096606} prove this conjecture, and they conjecture the  stronger statement: if the average is less than $\frac{16}{5}$, then $G$ is solvable. Moret\'o and Hung settle this in Theorem~A of \cite{MR3210700} by proving the conjecture. In addition, they show that $A_{5}$ meets this bound; so the bound is best. 

We now establish some more notation. Let $p$ be a prime, and take ${\rm Irr}_{p'}(G)$ to be the set of (complex) irreducible characters whose degrees are not divisible by $p$. Write ${\rm acd}_{p'}(G)$ for the average of the degrees of characters in ${\rm Irr}_{p'}(G)$. Similarly, when $k$ is a field, write ${\rm Irr}_{k}(G)$ for the set of irreducible characters in $G$ whose values lie in $k$ and ${\rm Irr}_{k,p'}(G)$ for the set of irreducible characters in $G$ with values lying in $k$ and whose degrees are not divisible by $p$. We write ${\rm acd}_{k}(G)$ for the average of the degrees of the characters in  ${\rm Irr}_{k}(G)$ and ${\rm acd}_{k,p'}(G)$ for the average of the degrees of the characters in ${\rm Irr}_{k,p'}(G)$. When $k=\mathbb{C}$, we use ${\rm acd}(G)$ in place of $ {\rm acd}_{\mathbb{C}}(G)$ and ${\rm acd}_{p'}(G)$ in place of ${\rm acd}_{\mathbb{C},p'}(G)$. 

Moret\'o and Hung look at this problem from another angle. They suggest in \cite{MR3210700} that one can look at a subset of these irreducible characters and still get information about the structure of a group. In particular, they define ${\rm acd}_{even}(G)$ to be the average of the degrees of the irreducible characters whose degrees are even. They show that if ${\rm acd}_{even}(G)<\frac{16}{5}$, then $G$ is solvable, see \cite{MR3210700}. Further, they establish the idea of looking at the irreducible characters having degrees that are not divisible by $p$ where $p$ is some prime that divides the order of $G$. They ask in that paper if the condition that ${\rm acd}_{3'}(G)<3$ or ${\rm acd}_{p'}(G)<\frac{16}{5}$ when $p>5$ implies that $G$ must be solvable. In Theorem 1.2 of \cite{MR3615444}, Hung settles this problem by showing: if ${\rm acd}_{3'}(G)<3$, ${\rm acd}_{5'}(G)<\frac{11}{4}$, or ${\rm acd}_{p'}(G)<\frac{16}{5}$ when $p>5$, then $G$ is solvable. He also presents examples that illustrate that these bounds are best possible. 

We study a condition on the average value of $k$-valued irreducible characters of a group of odd order that  guarantees that the group contains a normal $p$-complement for some prime $p$. A theorem of Thompson (Theorem 1 of \cite{MR166261} or Corollary 12.2 of \cite{MR1280461}) asserts that $G$ has a normal $p$-complement if the degree of every irreducible character of $G$ is 1 or divisible by $p$. This theorem is reformulated by Hung in \cite{MR3615444} as follows: if $G$ is a group and ${\rm acd}_{p'}(G)=1$, then $G$ has a normal $p$-complement. Furthermore, Hung improves Thompson's theorem in Theorem~1.1 of \cite{MR3615444}: if ${\rm acd}_{2'}(G) < 3/2$, then $G$ is $2$-nilpotent, and if ${\rm acd}_{p'}(G) < 4/3$, then $G$ is $p$-nilpotent when $p$ is an odd prime. Furthermore, Lewis pushes this further in Corollary 1.2 of \cite{MR3552294}. Let $G$ be a group and let $p$ be an odd prime. If ${\rm acd}_{p'}(G)<2(p+1)/(p+3)$, then $G$ is $p$-nilpotent.  

Let $\mathbb{Q}_p$ be the cyclotomic extension of $\mathbb{Q}$ by a $p^{th}$ root of unity. In Theorem 1.5 of \cite{MR3552294}, Lewis proves that if $G$ is a group of odd order and if $p=7$ and ${\rm acd}_{7'}(G)<9/5 $ or $p\ne7$ and ${\rm acd}_{p'}(G)<2 $ or ${\rm acd}_{\mathbb{Q}_{p},p'}(G)<2 $ or ${\rm acd}_{\mathbb{Q}_{p}}(G)<2 $, then $G$ is $p$-nilpotent. In this paper, we are able to find an upper bound that depends on the odd prime $p$ for any extension of $\mathbb{Q}$ that contains the primitive $p^{th}$ roots of unity. In particular,  for the fields $k=\mathbb{C}$ and $k=\mathbb{Q}_p$ we have the following: 
\begin{theorem} \label{thm1.1}
	Let $G$ be a group of odd order and $p$ an odd prime. Assume one of the following holds:	
	\begin{enumerate}
		
		\item ${\rm acd}(G)<3(p+2)/(p+8)~$ if $~p\equiv1\pmod3$; 
		\item  ${\rm acd}(G)<3(p^2+2)/(p^2+8)~$ if $~p\equiv2\pmod3$ and $(p-1)/2$ is even;
		\item ${\rm acd}(G)<3(p-1)/(p+3)$~~if~~$p\equiv2\pmod3$ ~~and ~~$(p-1)/2$ is odd;
		\item ${\rm acd}_{\mathbb{Q}_p}(G)<3p^2/(p^2+2)~$ if $~p\equiv2\pmod3$ and $(p-1)/2$ is odd;		
		\item ${\rm acd}_{\mathbb{Q}_p}(G)<3p/(p+2)~$ if $~p\equiv1\pmod3$; 
		\item  ${\rm acd}_{\mathbb{Q}_p}(G)<3p^2/(p^2+2)~$ if $~p\equiv2\pmod3$ and $(p-1)/2$ is even.
	\end{enumerate}
	Then $G$ is $p$-nilpotent.

\end{theorem}

This paper is organized as follows: in Section 1, we introduce the motivation and results related to the average character degree problem. Within Section 2, we state various basic results. In Sections 3 and 4, we prove the result in the critical restricted cases. Section 5 features the statement of the main result of the paper, the proof of the main theorem, and several examples.

\section{Basic Lemmas}
In this chapter, we present some basic results. The proof of Lemma~\ref{lem2.1} follows using the usual techniques from calculus for finding the minimum value of a continuous function on a closed interval. 

\begin{lemma}\label{lem2.1}
	Let $p$ be an odd prime number, let $a$ be a positive integer, and let $\displaystyle f(x) = \frac{x(x+p^a-1)}{x^2 + p^a -1}$, and $\displaystyle g(x) = \frac{xp^a}{x + p^a -1}$.
	\begin{enumerate}
		\item If $p^a\ge7$, then the minimum value for f(x) on $[3,\frac {p^a-1}2]$ is $\displaystyle \frac {3(p^a+2)}{p^a +8}$.
		\item If $p^a\ge11$, then the minimum value for $f (x)$ on $[5,\frac{p^a-1}{2}]$ is $\displaystyle \frac {3(p^a-1)}{p^a+3}$.
		\item If $p^a \ge 29$, then the minimum value for $f (x)$ on $[5,\frac{p^a-1}4]$ is $\displaystyle \frac {5(p^a+4)}{p^a+24}$.
		\item If $p^a\geq 11$, then the minimum value for $g(x)$ on $[3,\frac{p^a-1}{2}]$ is $\displaystyle \frac{3p^a}{p^a+2}$.
		\item If $p^a\geq 11$, then the minimum value for $g(x)$ on $[5,\frac{p^a-1}{2}]$ is $\displaystyle \frac{5p^a}{p^a+4}$.	
	\end{enumerate}

\end{lemma}

We need  the following result that we can prove by the same techniques as those used  in the proof of Lemma 2.3 of \cite{MR3552294}. 

\begin{lemma} \label{lem2.2}
	Let $p$ be an odd prime. Then the functions $f(x)=3(p^x+2)/(p^x+8)$, $g(x)=3p^x/(p^x+2)$  and $h(x)=3(p^x-1)/(p^x+3)$  are increasing functions on the interval $ [1,\infty)$ and so their minimum values on that interval are $f(1)=3(p+2)/(p+8)$, $g(1)=3p/(p+2)$, and $h(1)=3(p-1)/(p+3)$ respectively.
	

\end{lemma}

We now introduce some observations that appeared in Section 2 of \cite{MR3552294}. Let $k$ be a subfield of the complex numbers, and let $G$ be a group. Use $A^k(G)$ to denote the intersection of the kernels of the linear characters of $G$ with values in $k$, and note that  the set of linear characters with values in $k$ equals ${\rm Irr}(G/ A^{k}(G))$.  Furthermore, we can see that $G'\le A^k(G)$,  ${\rm Irr}_{k}(G)\cap {\rm Irr}(G/G')={\rm Irr}(G/A^{k}(G))$, and ${\rm Irr}_{k,p'}(G) \cap {\rm Irr}(G/G')={\rm Irr}(G/A^{k}(G))$. When $p$ is a prime, we define $A^p(G)$ to be the smallest normal subgroup of $G$ whose quotient is an elementary abelian $p$-group. Observe that by \cite{MR3552294} if $k=\mathbb{C}$, $\mathbb{Q}$, or $\mathbb{Q}_{p}$ where $p$ an odd prime, then  $A^k(G)=G'$, $A^k(G)=A^2(G)$, or  $A^k(G)=A^2(G) \cap A^p(G)$ respectively. 

We need the following results that have been proved by Lewis in Lemma 2.4 of \cite{MR3552294}. 

\begin{lemma} \label{lem2.3}
Let $K$ be a normal subgroup of $G$ so that $K\cap G'=1$. If $\varphi \in {\rm Irr}(K)$, then $\varphi$ extends to $\tilde{\varphi} \in {\rm Irr}(G)$ and the map $f: {\rm Irr}(G/K)\longrightarrow {\rm Irr}(G|\varphi)$ defined by $f(\chi)=\chi\tilde{\varphi}$ is a bijection such that $f(\chi)(1)=\chi(1)$ for all $ \chi \in {\rm Irr}(G/K)$.

\end{lemma}

We can use the same idea in the proof of Lemma 2.5 in \cite{MR3552294} and apply it for the groups of odd order to see that Lemma~\ref{lem2.4} is true. 

\begin{lemma} \label{lem2.4}
Suppose $k$ is a field, $G$ is a group of odd order, and $p$ is a prime. Assume $K$ is a minimal normal subgroup of $G$ such that $G'\cap K=1$. If ${\rm acd}_{k,p'}(G) \le3$, then ${\rm acd}_{k,p'}(G/K) \le{\rm acd}_{k,p'}(G)$.

\end{lemma}   
\section{When $H$ is abelian}
In Sections 3 and 4, we consider a group $G$ with odd order and having the following structure.  Suppose that $G = HV$ and $H \cap V = 1$ for subgroups $H$ and $V$ where $V$ is an elementary abelian normal $p$-group for some prime $p$ and $V$ can be viewed as a module for $H$ that is faithful and irreducible. In Lemma~\ref{lem3.1} of this section we will consider the case when $H$ is a nontrivial abelian group, and we will see that the lower bound for the average value of the degrees of the $k$-valued irreducible characters is $3(p+2)/(p+8)$ when $p\equiv1\pmod3$. We derive lower bounds in other circumstances. 

In the following lemma, we use the fact that when $H$ is abelian that the group $G$ is a Frobenius group. Observe that all of the values of the characters in $ {\rm Irr}(V)$ are $p^{th}$ roots of unity, so they have values in $k$. By Theorem 6.34 of \cite{MR1280461} the nonlinear irreducible characters of $G$ are induced from the non-principal irreducible characters of $V$. This implies that $G$ has exactly $(|V|-1)/|H|$ nonlinear $k$-valued irreducible characters and these all have degree equal to $|H|$ and vanish on $G\setminus V$ (see Problem 2.18 of \cite{MR1280461}). Furthermore, $|H|$ must divide $|V|-1$, which implies that $|H|$ is co-prime to $p$. By Theorem 6.15 of \cite{MR1280461}, this implies that every character in ${\rm Irr}_{k}(G)$ has $p'$-degree. In particular, we can write that ${\rm Irr}_{k,p'}(G)={\rm Irr}_{k}(G)$ and ${\rm acd}_{k,p'}(G)={\rm acd}_{k}(G)$. 

Suppose $|V|=p^a = 3$ or $p^a = 5$. We know $|H|$ divides $|V|-1 = p^a -1$, and $p^a -1$ is a power of 2 in these cases. On the other hand, $|H|$ is odd and not 1. Therefore, this leads to a contradiction.  Hence, we may assume that $p^a \ge 7$.

Note that we set $t$ to be the least common multiple of the orders of the roots of unity in $k$ having order dividing $|G|$. We see that the value ${\rm acd}_{k}(G)$ has different lower bounds depending on the odd prime $p$. 
	
\begin{lemma}\label{lem3.1}
	Assume $H$ acts faithfully on an irreducible module $V$ of characteristic $p$, and suppose that $G=HV$ is a group of odd order. Assume that $k$ is an extension of $\mathbb{Q}$ that contains the primitive $p^{th}$ roots of unity. Suppose that $t$ is the least common multiple of orders of roots of unity in $k$ having order dividing $|H|$. Let H be a nontrivial abelian group. Then

		\begin{enumerate}
			\item ${\rm acd}_{k}(G)\ge\frac{3(p+2)}{p+8}~~~~~~$ if $~~~~~p\equiv1\pmod{3}$ and $~3$ divides $t$,
			\item ${\rm acd}_{k}(G)\ge\frac{3(p^2+2)}{p^2+8}~~~~~$ if $~~~~~p\equiv2\pmod{3}$, $~\frac{p-1}{2}\nmid t~$ and $~3$ divides $t$, 
			\item ${\rm acd}_{k}(G)\ge\frac{3(p-1)}{p+3}~~~~~~$ if $~~~~~p\equiv2\pmod{3}$ and $~\frac{p-1}{2}\mid t$,
			\item ${\rm acd}_{k}(G)\ge\frac{3p}{p+2}~~~~~~~~~$ if $~~~~~p\equiv1\pmod{3}$ and $3$ does not divide $t$,
			\item ${\rm acd}_{k}(G)\ge\frac{3p^2}{p^2+2}~~~~~~~~$ if $~~~~~p\equiv2\pmod{3}$, $~\frac{p-1}{2}\nmid t~$ and $3$ does not divide $t$. 
	\end{enumerate}
	
	\begin{proof}
		We have $G/V\cong H$. Since $H$ is abelian, this implies that $G'\le V$. Because $H$ acts faithfully on $V$, we know that $G'\neq1$, and because $V$ is irreducible under the action of $H$, we have $G'=V$. Moreover, we see that $H$ is cyclic by Lemma 0.5 of \cite{MR1261638}. Hence, $G$ is a Frobenius group in this case. Now, Lemma 3.1 of \cite{MR3552294} gives us  the average value of the degrees of the $k$-valued irreducible characters by the following equation $${\rm acd}_{k}(G) =\frac{|H|(|H:A^k(H)|+p^a-1)}{|H||H:A^k(H)|+p^a-1}.$$ 
		
		Suppose that $3$ divides $t$. Set $l = |H:A^k (H)|$.  Now, because $H$ is cyclic, there is a linear character of $H/A^k(H)$ whose order is $l$, and it is clear that this character must have at least one value that is a root of unity of order $l$.  It follows that $l$ will divide both $|H|$ and $t$.  Let $d=gcd(|H|,t)$, and we have just shown that $l$ divides $d$.  Since $H$ is cyclic, it follows that $H$ has a linear character $\lambda$ of order $d$.  It is not difficult to see that the values of $\lambda$ will all be roots of unity of orders dividing $d$.  Since $d$ divides $t$, we deduce that they all will lie in $k$. Hence, the values of $\lambda$ lie in $k$ and so, $A^k(H) \le \ker(\lambda)$, and so, $d$ divides $l$.  We conclude that $l=d$.
		
		We know that $|H:A^k(H)|$ is the number of the linear characters in $H$ that have values in $k$. This implies that the index $|H:A^k(H)|$ is the greatest number that divides $|H|$ and the number of the $t^{th}$ roots of unity in $k$. Then we can say that $|H:A^k(H)|$ is the greatest common divisor of $t$ and $|H|$. Hence $l=gcd(t,|H|)$. 
		
		The equation of ${\rm acd}_k(G)$ becomes $ {\rm acd}_{k}(G)=|H|(l+p^a-1)/(l|H|+p^a-1)$, where $1\le l\le |H|\le(p^a-1)/2 $. For a fixed $l$, the function $g(x)=x(l+p^a-1)/(lx+p^a-1)$ is increasing. Assume that $l=1$, so that $t$ is co-prime to $|H|$. Since 3 divides $t$, we conclude that $3\nmid|H|$, and hence, $|H|\ge5$ and $p^a\ge 11$. We see that $g(x)\ge 5p^a/(p^a+4)$ by Lemma~\ref{lem2.1}(5). Also, if $l\ne1$, then we have $g(x)\ge l(l+p^a-1)/(l^2+p^a-1)$. Hence, we obtain
		
		$$\displaystyle {\rm acd}_{k}(G) \ge \left\{
		\begin{array}{ll}
			\frac{l(l+p^a-1)}{l^2+p^a-1} ~~~~~ if ~~~~~ 3\le l\le \frac{p^a-1}{2} \\\\
			
			\frac{5p^a}{p^a+4} ~~~~~~~~~ if ~~~~~~~~~ l=1 \\
		\end{array}
		\right. $$
		
		We define the function $f$ in the single variable $l$ as follows: $f(l)=l(l+p^a-1)/(l^2+p^a-1)$. We consider two cases depending on whether $p$ is congruent to $1$ or $2\pmod{3}$.  The first case is that $p\equiv1\pmod{3}$. We now consider $l\ge3$. We see that ${\rm acd}_k(G)\ge3(p^a+2)/(p^a+8)$ by Lemma~\ref{lem2.1}(1). Second, suppose that $l=1$. We have that ${\rm acd}_{k}(G)\ge5p^a/(p^a+4)$. Now, for any values of $p^a$, this yields $3(p^a+2)/(p^a+8)\le5p^a/(p^a+4)$. Using Lemma~\ref{lem2.2}, we have in this case the desired result that ${\rm acd}_{k}(G)\ge3(p+2)/(p+8)$. 
		
		Now, assume $p\equiv2\pmod{3}$. We have seen that the case $p^a=5$ does not occur. Therefore, we must have $p^a\ge11$. First, we handle the case $l=3$. Substituting into the above function we obtain, $f(l)=3(p^a+2)/(p^a+8)$. Since $3$ does not divide $p-1$, it follows that $|V|\ge p^2$.  This implies that ${\rm acd}_{k}(G)\ge3(p^2+2)/(p^2+8)$ by Lemma~\ref{lem2.2} for the interval $ [2,\infty)$.
		
		Suppose that $(p-1)/2$ does not divide $t$. Thus, $(p-1)/2$ does not divide $l$ and we have two cases for $l$. First, consider $l=1$. Therefore, we have that ${\rm acd}_{k}(G)\ge5p^a/(p^a+4)$. By Lemma~\ref{lem2.2}, ${\rm acd}_{k}(G)\ge5p/(p+4)$, and since $3(p^2+2)/(p^2+8)\le5p/(p+4)$, we have ${\rm acd}_{k}(G)\ge3(p^2+2)/(p^2+8)$. Second, consider $5\le l< (p-1)/2$. Recall that $l$ divides $p-1$ and 3 does not divide $p-1$, so $l \le (p-1)/4$. This implies that $5\le l\le(p-1)/4$. By Lemma~\ref{lem2.1}(3), we compute that $ {\rm acd}_{k}(G)\ge5(p+4)/(p+24)$. Again, for any value $p$ we have $3(p^2+2)/(p^2+8)<5(p+4)/(p+24)$. Hence, we deduce that ${\rm acd}_{k}(G)\ge3(p^2+2)/(p^2+8)$.
		
		Suppose that $(p-1)/2$ divides $t$. Then $(p-1)/2$ is odd. Thus, we have only the case that $5\le l\le(p-1)/2$. By Lemma~\ref{lem2.1}(2) and Lemma~\ref{lem2.2}, we can calculate that $f(l)\ge3(p-1)/(p+3)$. Thus, we have ${\rm acd}_{k}(G)\ge3(p-1)/(p+3)$. Since $3(p-1)/(p+3)<3(p^2+2)/(p^2+8)$, we conclude that ${\rm acd}_{k}(G)\ge3(p-1)/(p+3)$. 
		
		For the rest of this proof, assume that 3 does not divide $t$.  We still 
		use the notation $|H:A^k (H)|=l$ where $l$ is the greatest common divisor of $|H|$ and $t$. Note since $3$ does not divide $t$, that $l$ is not divisible by $3$. As above, we have that ${\rm acd}_{k}(G)=|H|(l+p^a-1)/(l|H|+p^a-1)$ and $\max(3,l)\le |H|\le(p^a-1)/2$.
		
		For a fixed $l$, the function $g(x)=x(l+p^a-1)(lx+p^a-1)$ is increasing on $[3,(p^a-1)/2]$. Consider that $l=1$. This implies that $3\le|H|\le(p^a-1)/2$. Thus, the minimum value for $g(x)$ on the interval $[3,(p^a-1)/2]$ is $3p^a/(p^a+2)$ by Lemma~\ref{lem2.1}(4). Now, consider that $l\ne1$, so $l\le|H|\le(p^a-1)/2$, and the minimum value for $g(x)$ on the interval $[l,(p^a-1)/2]$ is $l(l+p^a-1)/(l^2+p^a-1)$. Hence, we can write that 
		$${\rm acd}_{k}(G)\ge \left\{
		\begin{array}{ll}
			\frac{l(l+p^a-1)}{l^2+p^a-1}  ~~~~~~~ if ~~~~~~ 5\le l\le\frac{p^a-1}{2} \\\\
			\frac{3p^a}{p^a+2} ~~~~~~~~~~~ if ~~~~~~~~ l=1 \\
		\end{array}
		\right. $$
		
		Again, we need to find the best lower bound in this case. Define the function $f$ in single variable $l$ as follows: $f(l)=l(l+p^a-1)/(l^2+p^a-1)$. First, assume $p\equiv1\pmod{3}$. Since $l$ is not divisible by $3$, we observe that $(p-1)/2$ does not divide $l$. This implies that we have two cases for $l$. First, consider $l=1$. We compute that ${\rm acd}_{k}(G)\ge3p/(p+2)$ by Lemma~\ref{lem2.2}. Now, suppose that $5\le l<(p-1)/2$. Thus, we consider that $5\le l\le(p-1)/4$. Now, by Lemma~\ref{lem2.1}(3) and by Lemma~\ref{lem2.2} we see that $f(l)\ge5(p+4)/(p+24)$. This yields ${\rm acd}_{k}(G)\ge5(p+4)/(p+24)$. And, since $3p/(p+2)\le5(p+4)/(p+24)$ for all values of $p$, we obtain the desired result that ${\rm acd}_{k}(G)\ge3p/(p+2)$. 
		
		Finally, assume that $p\equiv2\pmod{3}$. We first handle the case when $l=1$. Thus, ${\rm acd}_k(G)\ge3p^a/(p^a+2)$. Since $G$ is a Frobenius group and $3$ does not divides $p-1$, we know that $|V|\ge p^2$. By Lemma~\ref{lem2.2}, we can see that ${\rm acd}_k(G)\ge3p^2/(p^2+2)$.
		
		We now suppose $5\le l$.  There are two cases. The first case is $(p-1)/2$ does not divide $t$. As before, in this case, we have $l\le(p-1)/4$. By Lemma~\ref{lem2.1}(3) and Lemma~\ref{lem2.2}, we see that $f(l)\ge5(p+4)/(p+24)$. Since $3p^2/(p^2+2)\le5(p+4)/(p+24)$, we may conclude that ${\rm acd}_k(G)\ge3p^2/(p^2+2)$.
		
		The second case is $(p-1)/2$ does divide $t$.  In this case, we have $l\le(p-1)/2$. Thus, Lemma~\ref{lem2.1}(2) and Lemma~\ref{lem2.2} give us that $f(l)\ge3(p-1)/(p+3)$. Since $3(p-1)/(p+3)\le3p^2/(p^2+2)$, we conclude that ${\rm acd}_k(G)\ge3(p-1)/(p+3)$.	
	\end{proof}	
\end{lemma}

%
%

We now consider the odd prime $p=3$ in the following lemma and we will use a different method to prove our results.  

\begin{lemma} \label{lem3.2}
	Assume $H$ acts faithfully on an irreducible module $V$ of characteristic $3$, and suppose that $G=HV$ is a group of odd order. Suppose that $k$ is an extension of $\mathbb{Q}$ that contains the primitive $3^{rd}$ roots of unity. Let $H$ be a nontrivial abelian group. 
	Then, ${\rm acd}_k(G)\ge\sfrac{182}{61}$ except if  $A^k(H)=1$  and $|H|=13$, then ${\rm acd}_k(G)\ge\sfrac{13}{5}$.
	
	
	\begin{proof}
		Recall the equation of the average value of the degrees of the $k$-valued irreducible characters and set $p=3$: $${\rm acd}_{k}(G)=\frac{|H|(|H:A^k(H)|+3^a-1)}{|H||H:A^k(H)|+3^a-1} .$$
		
		Now, notice that $H$ acts faithfully and irreducibly on $V$, and $|H|$ divides $|V|-1$. Since $|V|=3^a$, this implies that $|H|$ is not divisible by $3$, and so, $5\le|H|\le(3^a-1)/2$. First, we handle the case when $A^k(H)=1$. In this case the equation becomes ${\rm acd}_k(G)=|H|(|H|+3^a-1)/(|H|^2+3^a-1)$ and the function $f(x)=x(x+3^a-1)/(x^2+3^a-1)$ is increasing on the interval $5\le x\le (3^a-1)/2$. 
		Consider that $|H|=5$. This implies $|V|=3^a\ge3^4$ and $f(5)=5(3^a+4)/(3^a+24)$ is at least $5(3^4+4)/(3^4+24)=85/21$ for all $a\ge 4$. We can conclude that  ${\rm acd}_k(G)\ge 85/21\ge 182/61$. Consider $|H|=7$ or $11$. This implies $|V|\ge3^6$ or $3^5$. In the same way we can see that ${\rm acd}_k(G)\ge81/17\ge182/61$ in both cases. Now, assume $|H|=13$. This implies $|V|\ge3^3$. We can show that ${\rm acd}_k(G)\ge13/5$ which is the exception. Finally, assume $|H|\ge17$. This implies $|V|\ge3^{6}$  and the function $f(x)=x(x+3^a-1)/(x^2+3^a-1)$ has a minimum value $3(3^a-1)/(3^a+3)$ at $x=(3^a-1)/2$. Again, we see that ${\rm acd}_k(G)\ge182/61$.
		
		
		
		For the rest of the proof, assume that $|A^k(H)|\ne1$. Notice, since $|H|$ is not divisible by $3$, that we have $|A^k(H)|\ge5$. Consider $|H|=5$ or $7$. This implies ${\rm acd}_k(G)=|H|3^a/(|H|+3^a-1)$. Then, we see that ${\rm acd}_k(G)=5(3^a)/(3^a+4)$ or ${\rm acd}_k(G)=7(3^a)/(3^a+6)$. So, we obtain the desired result ${\rm acd}_k(G)\ge182/61$ in both cases. Now, assume that $|H|\ge11$~(we may assume $|H|\ne9$, because $|H|$ is not divisible by 3). We can write $(187/81)|H|-11\ge106/81|H|$, and $(|H|-81/17)\ge106/187|H|$. Also, we have that $|A^k(H)|\ge5$ and $(3^a-1)/|H|\ge2$.  This implies $(3^a-1)\ge10|H:A^k(H)|$. Thus, we compute $$\displaystyle (|H|-81/17)(3^a-1)\ge\sfrac{64}{17}|H||H:A^k(H)| .$$
		
		Finally, we compute \begin{align*} |H|(|H:A^k(H)|+3^a-1)&=|H||H:A^k(H)|+(|H|-\sfrac{81}{17})(3^a-1)+\sfrac{81}{17}(3^a-1)\\ &\ge |H||H:A^k(H)|+\sfrac{64}{17}|H||H:A^k(H)|+\sfrac{81}{17}(3^a-1)\\ &=\sfrac{81}{17}(|H||H:A^k(H)|+3^a-1). 
		\end{align*} 
		
		Hence, we obtain the desired result. 
	\end{proof}
\end{lemma}

\section{When $H$ is nonabelian}
In Section 4, we present some lemmas to handle the case when $H$ is a nonabelian subgroup of the group $G=HV$ that appeared early in Section 3. In the following lemma we will see that every $H$-orbit on ${\rm Irr}(V)\setminus\{1_V\}$ has length at least $7$, and at least one such orbit has length that is not divisible by $p$. 

\begin{lemma}\label{lem4.1}
	Assume $H$ acts faithfully on an irreducible module $V$ of characteristic $p$, and suppose that $G=HV$ is a group of odd order. Let $H$ be a nonabelian group. Then, every $H$-orbit on ${\rm Irr}(V)\setminus\{1_V\}$ has length at least $7$, and at least one such orbit has length that is not divisible by $p$. Furthermore, among the $H$-orbits whose length is not divisible by $p$, either there are at least two $H$-orbits of length $7$ or there is at least one $H$-orbit of length at least $9$, and when $p=3$, this orbit must have length at least 11.  
	
	\begin{proof}
		The set ${\rm Irr}(V)$ may be viewed as a faithful irreducible module for $H$ of order $p^{a}$. Suppose $\alpha \in{\rm Irr}(V) \setminus\{1_V\}$, and let $T$ be the stabilizer of $\alpha$ in $H$. Since $H$ acts faithfully and irreducibly on $V$, we see that $H$ acts faithfully and irreducibly on ${\rm Irr}(V)$. In particular, ${\rm Irr}(V)$ as a module will be generated by the elements in the orbit of $\alpha$. Since the action of $H$ is faithful, $\cap_{h\in H}T^h=1$ and thus, $T$ is core-free in $H$. Also, we see that $H$ is isomorphic to a subgroup of $S_{n}$ where $n=|H:T|$ and $n$ is odd. Since $H$ is nonabelian, it must be that $n>2$. Also, $n\ne3$ or $5$, because $S_3$ and $S_5$ do not have any nonabelian subgroups of odd order. This implies that $n=|H:T|$ is at least $7$. Hence, every $H$-orbit on ${\rm Irr}(V)\setminus\{1_V\}$ has length at least $7$. Also, the lengths of the $H$-orbits of ${\rm Irr}(V)\setminus\{1_V\}$ must sum to $p^{a}-1$, so that at least one of these orbits must have length that is not divisible by $p$. 
		
		Assume that the last conclusion is not true. This implies that there is only one $H$-orbit whose length is not divisible by $p$ and the length of this orbit is equal $7$. Let $\alpha$ be a representative of the $H$-orbit whose length is not divisible by $p$. Let $T$ be the stabilizer of $\alpha$ in $H$ and we know that $|H:T|=7$. Now, we know that ${\rm Irr}(V)$ is the union of $H$-orbits, the orbit of the principal character $1_V$, the orbit of $\alpha$, and the remaining orbits (whose lengths are divisible by $p$). This yields $p^a=|V|=1+7+pr$ for some positive integer $r$, and then $p^a\equiv8\pmod p$. This implies that $p$ divides $8$ which is a contradiction since $p$ is an odd prime. Note that the smallest length of an orbit whose length is not divisible by 3 will not be 9, so it must be that 11 is the smallest such length when $p=3$.
	\end{proof}	
\end{lemma}
In Theorem 3.3 of \cite{MR3552294}, Lewis demonstrates conclusion (1) of the following lemma, and the proofs of conclusions (2) and (3) are also motivated by the proof of Theorem 3.3 of \cite{MR3552294}. For any prime $p$ and a field $k$, define ${\rm nl}_{k,p'}(G)$ to be the set of nonlinear irreducible characters in $G$ with values in $k$ and whose degree are not divisible by $p$. 
\begin{lemma}\label{lem4.2}
	Assume $H$ acts faithfully on an irreducible module $V$ of characteristic $p$, and suppose that $G=HV$ is a group of odd order. Suppose that $k$ is an extension of $\mathbb{Q}$ that contains the primitive $p^{th}$ roots of unity. Let $H$ be a nonabelian group. Then
	\begin{enumerate}
		\item $\displaystyle|{\rm Irr}_{k,p'}(G)| = |H:A^k(H)|+|{\rm nl}_{k,p'}(H)|+\sum_{i=1}^{l}|{\rm Irr}_{k,p'}(T_{i})|$,
		\item $\displaystyle  \sum_{\chi \in {\rm Irr}_{k,p'}(G)}\chi(1) \geq |H:A^k(H)| + 3|{\rm nl}_{k,p'}(H)| + \sum_{i=1}^{l}|H:T_{i}||{\rm Irr}_{k,p'}(T_{i})|$, $p\ne3$,
		\item $\displaystyle \sum_{\chi \in {\rm Irr}_{k,3'}(G)}\chi(1) \geq |H:A^k(H)| + 5|{\rm nl}_{k,3'}(H)| + \sum_{i=1}^{l}|H:T_{i}||{\rm Irr}_{k,3'}(T_{i})|$.
		
	\end{enumerate}
	
	\begin{proof}
		By Lemma~\ref{lem4.1}, $H$ has an orbit on ${\rm Irr}(V) \setminus\{1_V\}$ whose length is not divisible by $p$ and that length is at least $7$. Let $\alpha_{1}, \alpha_{2},...,\alpha_{l}$ be representatives for the $H$-orbits of ${\rm Irr}(V) \setminus\{1_V\}$ whose lengths are not divisible by $p$. Let $T_{i}$ be the stabilizer of $\alpha_{i}$ in $H$. Furthermore, we apply Theorem 1.4 of \cite{MR2426855} (the Fundamental Counting Principle) to conclude that $|H:T_{i}|$ will equal the length of the $H$-orbit of $\alpha_{i}$. 
		
		
		Now arguing as in Theorem~3.3 of \cite{MR3552294}, we obtain Conclusion (1): \[|{\rm Irr}_{k,p'}(G)|=|H:A^k(H)|+|{\rm nl}_{k,p'}(H)|+\sum_{i=1}^{l}|{\rm Irr}_{k,p'}(T_{i})|.\]
		
		Using Theorem 6.11 of \cite{MR1280461}, we have $ {\rm Irr}_{k,p'}(G|\alpha_{i}) = \{(\alpha_{i} \times \tau)^{G} | \tau \in {\rm Irr}_{k,p'}(T_{i})\} .$ Since  $|G:VT_{i}| = |H:T_{i}|$, it follows that for any  $ \psi \in {\rm Irr}_{k,p'}(G|\alpha_{i})$  we have $\psi(1) = |H:T_{i}|\tau(1)$ where $\tau \in {\rm Irr}_{k,p'}(T_{i})$. Thus, we obtain $$\displaystyle \sum_{\chi \in {\rm Irr}_{k,p'}(G)}\chi(1) = \sum_{\gamma \in {\rm Irr}_{k,p'}(H)}\gamma(1)+\sum_{i=1}^{l}|H:T_{i}|(\sum_{\theta \in {\rm Irr}_{k,p'}(T_i)}\theta(1)).$$
		Also, since $|H|$ is odd and $\gamma(1)$ divides $|H|$, we see that $\gamma(1)\ge3$ for all $\gamma \in {\rm nl}_{k,p'}(H)$. Thus, $$ |H:A^k(H)| + \sum_{\gamma(1)\neq 1, \gamma \in {\rm Irr}_{k,p'}(H)}\gamma(1) = \sum_{\gamma \in {\rm Irr}_{k,p'}(H)}\gamma(1) \geq |H:A^k(H)|+ 3|{\rm nl}_{k,p'}(H)| ,$$ and $$ \sum_{\theta \in {\rm Irr}_{k,p'}(T_{i})}\theta(1) \geq |{\rm Irr}_{k,p'}(T_{i})| .$$
		This yields Conclusion (2): $$ \sum_{\chi \in {\rm Irr}_{k,p'}(G)}\chi(1) \geq |H:A^k(H)| + 3|{\rm nl}_{k,p'}(H)| + \sum_{i=1}^{l}|H:T_{i}||{\rm Irr}_{k,p'}(T_{i})|.$$
		
		Finally, assume $p=3$. From the above work, it is still true that $$|{\rm Irr}_{k,3'}(G)|=|H:A^k(H)|+|{\rm nl}_{k,3'}(H)|+\sum_{i=1}^{l}|{\rm Irr}_{k,3'}(T_{i})|,$$ and $$\displaystyle \sum_{\chi \in {\rm Irr}_{k,3'}(G)}\chi(1) = \sum_{\gamma \in {\rm Irr}_{k,3'}(H)}\gamma(1)+\sum_{i=1}^{l}|H:T_{i}|(\sum_{\theta \in {\rm Irr}_{k,3'}(T_i)}\theta(1)).$$
		
		Since $H$ acts faithfully on the irreducible module $V$ of characteristic 3, then we have $|H|\ge5$, and $\gamma(1)$ divides $|H|$. This implies $\gamma(1)\ge5$ for all $\gamma \in {\rm nl}_{k,3'}(H)$. Thus, we compute that  $$ |H:A^k(H)| + \sum_{\gamma(1)\ne1, \gamma \in {\rm Irr}_{k,3'}(H)}\gamma(1) = \sum_{\gamma \in {\rm Irr}_{k,3'}(H)}\gamma(1) \ge|H:A^k(H)|+ 5|{\rm nl}_{k,3'}(H)| ,$$ 
		and observe that $$\sum_{\theta \in {\rm Irr}_{k,3'}(T_{i})}\theta(1) \ge|{\rm Irr}_{k,3'}(T_{i})| .$$ 
		
		Hence, Conclusion (3) follows: $$ \sum_{\chi \in {\rm Irr}_{k,3'}(G)}\chi(1) \ge|H:A^k(H)| + 5|{\rm nl}_{k,3'}(H)| + \sum_{i=1}^{l}|H:T_{i}||{\rm Irr}_{k,3'}(T_{i})|.$$				
	\end{proof}	
\end{lemma}
Now, we are ready to present the proof of the case where $H$ is nonabelian. In the following lemma we will consider the nonabelian case for the subgroup $H$. And we will show that the ${\rm acd}_{k,p'}(G)\ge3$ for $p\ne3$, and  ${\rm acd}_{k,3'}(G)\ge81/17$. These values are the best lower bounds that we found for any field $k$ such that $\mathbb{Q}_p \subseteq k\subseteq \mathbb{C}$. We can do better for other circumstances. 

\begin{lemma}\label{lem4.3}
	Assume $H$ acts faithfully on an irreducible module $V$ of characteristic $p$, and suppose that $G=HV$ is a group of odd order. Suppose that $k$ is an extension of $\mathbb{Q}$ that contains the primitive $p^{th}$ roots of unity. Let $H$ be a nonabelian group. Then the following holds. 
	\begin{enumerate}
		\item If $p\ne3$, then ${\rm acd}_{k,p'}(G)\geq3$.
		\item If $p=3$, then ${\rm acd}_{k,3'}(G)\ge\sfrac{81}{17}$.
	\end{enumerate}
	\begin{proof}
		Working as in Lemma~\ref{lem4.1}, take $\alpha_{1}, \alpha_{2},...,\alpha_{l}$ to be representatives for the $H$-orbits of ${\rm Irr}(V) \setminus\{1_V\}$ whose lengths are not divisible by $p$. For all $i$, let $T_{i}$ be the stabilizer of $\alpha_{i}$ in $H$  and note that  $|H:T_i|$ is not divisible by $p$ and $|H:T_i|\ge7$. Also, we have that $T_i$ is core-free in $H$ for all $i$. 
		
		First, we notice that $ A^k(T_1)<A^k(H) $. Otherwise, we have $A^k(T_1)=A^k(H)$, and this implies that $1\le A^k(H) \le T_1$ which contradicts the fact that $T_1$ is core-free. Thus, we have $|H:A^k(T_1)|\ge3|H:A^k(H)|$. Note that the number of linear characters of $T_1$ with values in $k$ is $|T_1:A^k (T_1)|$, and the linear characters of $T_1$ with values in $k$ all lie in ${\rm Irr}_{k,p'}(T_1)$. This yields $|{\rm Irr}_{k,p'} (T_1)| \ge |T_1:A^k (T_1)|$.
		
		Now, we will work by cases on the length of $H$-orbits on $V$ whose lengths are not divisible by $p$. As shown in Lemma~\ref{lem4.1}, we have only two cases for $p\ne3$:
		
		
		
		\textbf{Case(1):} There are at least two $H$-orbits whose lengths are not divisible by $p$ and their lengths equal $7$. Among the $\alpha_{i}$ relabel the subscripts so that $|H:T_i|=7$ and $|{\rm Irr}_{k,p'}(T_{1})|\le|{\rm Irr}_{k,p'}(T_{2})|$ for $i=1$, $2$. Thus, we have that $3|H:T_1|=7+2|H:T_1|$ and $\displaystyle |H:T_1|-\frac{7}{3}=\frac{2}{3}|H:T_1|$. We can combine this with the previous facts that $|{\rm Irr}_{k,p'}(T_{1})|\ge|T_{1}:A^k(T_1)|$ and $\frac{1}{3}|H:A^k(T_1)|\ge|H:A^k(H)|$ to calculate the following 
		$$ \displaystyle \left(|H:T_{1}|-\frac{7}{3}\right) |{\rm Irr}_{k,p'}(T_{1})| \ge \frac{2}{3}|H:A^k(T_1)|\ge 2|H:A^k(H)|. $$
		
		In addition we compute 
		\begin{align*}
			|H:T_{1}||{\rm Irr}_{k,p'}(T_{1})|
			&=\left(|H:T_{1}|-\frac{7}{3}\right)|{\rm Irr}_{k,p'}(T_{1})|+ \frac{7}{3}|{\rm Irr}_{k,p'}(T_{1})| \\
			&\ge 2|H:A^k(H)| + \frac{7}{3}|{\rm Irr}_{k,p'}(T_{1})|.
		\end{align*}
		
		The existence of $\alpha_2$ guarantees that $\displaystyle \sum_{i=2}^{l}|{\rm Irr}_{k,p'}(T_{i})|\ne0$. With this in mind and the fact that $|H:T_i|\ge7$ for all $i$, we can calculate
		
		\begin{align*} 
			\sum_{i=1}^{l}|H:T_{i}||{\rm Irr}_{k,p'}(T_{i})| 
			&\ge 2|H:A^k(H)| + \sfrac{7}{3}|{\rm Irr}_{k,p'}(T_{1})| + 7\sum_{i=2}^{l}|{\rm Irr}_{k,p'}(T_{i})|\\
			&\ge 2|H:A^k(H)| + 3\sum_{i=1}^{l}|{\rm Irr}_{k,p'}(T_{i})|.	 
		\end{align*}
		
		Finally, by Lemma~\ref{lem4.2}(2) we compute 
		\begin{align*}
			\displaystyle \sum_{\chi \in {\rm Irr}_{k,p'}(G)} \chi(1) 
			& \ge|H:A^k(H)|+3|{\rm nl}_{k,p'}(H)|+\sum_{i=1}^{l}|H:T_{i}||{\rm Irr}_{k,p'}(T_i)|\\
			& \ge|H:A^k(H)|+3|{\rm nl}_{k,p'}(H)|+2|H:A^k(H)|+3\sum_{i=1}^{l}|{\rm Irr}_{k,p'}(T_{i})|\\
			& \ge3\left(|H:A^k(H)|+|{\rm nl}_{k,p'}(H)|+\sum_{i=1}^{l}|{\rm Irr}_{k,p'}(T_{i})| \right).
		\end{align*}
		Hence, by applying Lemma~\ref{lem4.2}(1) we get the desired result in this case that ${\rm acd}_{k,p'}(G)\ge3$. \label{case(2)}
		
		\textbf{Case(2):} There is at least one $H$-orbit whose length is not divisible by $p$ and that length is at least $9$. Among the $\alpha_{i}$, relabel the subscripts so that $|H:T_1|\ge9$. Thus, we can write that $3|H:T_{1}|\ge9 + 2|H:T_{1}|$ and $\displaystyle |H:T_{1}|-3\ge \sfrac{2}{3}|H:T_{1}|$. Recall that $|{\rm Irr}_{k,p'}(T_{1})|\ge|T_{1}:A^k(T_1)|$ and $\sfrac{1}{3}|H:A^k(T_1)|\ge|H:A^k(H)|$. Thus, we derive $$ \displaystyle (|H:T_{1}|-3) |{\rm Irr}_{k,p'}(T_{1})|\ge \sfrac{2}{3}|H:A^k(T_1)|\ge 2|H:A^k(H)|. $$ 
		
		Also, we calculate  
		\begin{align*}
			|H:T_{1}||{\rm Irr}_{k,p'}(T_{1})|
			&=(|H:T_{1}|-3)|{\rm Irr}_{k,p'}(T_{1})|+ 3|{\rm Irr}_{k,p'}(T_{1})|\\
			&\ge 2|H:A^k(H)| + 3|{\rm Irr}_{k,p'}(T_{1})|.
		\end{align*} 
		
		Notice that for all $i$ we have $|H:T_i|\ge7$, and then we have
		\begin{align*}
			\sum_{i=1}^{l}|H:T_{i}||{\rm Irr}_{k,p'}(T_{i})|
			& \ge2|H:A^k(H)|+3|{\rm Irr}_{k,p'}(T_{1})|+7\sum_{i=2}^{l}|{\rm Irr}_{k,p'}(T_{i})|\\
			& \ge2|H:A^k(H)|+3\sum_{i=1}^{l}|{\rm Irr}_{k,p'}(T_{i})|.	 
		\end{align*}
		
		Thus, by Lemma~\ref{lem4.2}(2) we compute 
		\begin{align*}
			\displaystyle \sum_{\chi \in {\rm Irr}_{k,p'}(G)} \chi(1) 
			& \ge|H:A^k(H)|+3|{\rm nl}_{k,p'}(H)|+\sum_{i=1}^{l}|H:T_{i}||{\rm Irr}_{k,p'}(T_i)|\\
			& \ge|H:A^k(H)|+3|{\rm nl}_{k,p'}(H)|+2|H:A^k(H)|+3\sum_{i=1}^{l}|{\rm Irr}_{k,p'}(T_{i})|\\
			& \ge3\left(|H:A^k(H)|+|{\rm nl}_{k,p'}(H)|+\sum_{i=1}^{l}|{\rm Irr}_{k,p'}(T_{i})|\right).
		\end{align*}
		
		In a similar way we obtain the desired result ${\rm acd}_{k,p'}(G)\ge3.$  
		
		
		For the rest of the proof, assume that $p=3$. In similar way, we will work by cases on the length of such $H$-orbits on $V$ whose lengths are not divisible by $3$. Again, by Lemma~\ref{lem4.1} we have the following two cases.
		
		
		\textbf{Case(1):} Consider that we have at least two $H$-orbits whose lengths are not divisible by 3 and their lengths are equal to 7. Among the $\alpha_{i},$ relabel the subscripts so that $|H:T_i|=7$ for $i=1$, $2$. In the second paragraph of this proof, we showed that $A^k(T_i)<A^k(H)$. So by The Diamond Isomorphism Theorem of \cite{MR2286236} we can derive that $|A^k(H):A^k(T_i)|\ge7$ and $|H:A^k(T_i)|\ge7|H:A^k(H)|$ for $i=1$, $2$. Also, the fact that $|H:T_i|=7$ implies $|H:T_i|\ge\frac{7(81)}{87}$ for $i=1$, $2$. Then for $i=1$, $2$, we see that $\frac{87}{17}|H:T_i|-\frac{7(81)}{17}\ge0$ and $7|H:T_i|-\frac{7(81)}{17}\ge\frac{32}{17}|H:T_i|$. This implies the following inequality for $i=1$, $2$:  $$|H:T_i|-\frac{81}{17}\ge \frac{32}{7(17)}|H:T_i|.$$
		Also, for $i=1$, $2$, we know that every linear character of $T_i$ with values in $k$ is contained inside ${\rm Irr}_{k,3'}(T_i)$, so $|{\rm Irr}_{k,3'}(T_i)|\ge|T_i:A^k(T_i)|$. 
		
		Now, we combine these to observe that
		\begin{align*}
			\left(|H:T_i|-\frac{81}{17}\right)|{\rm Irr}_{k,3'}(T_i)|&\ge \frac{32}{7(17)}|H:T_i||T_i:A^k(T_i)|\\&=\frac{32}{7(17)}|H:A^k(T_i)|\\&\ge \frac{32}{17}|H:A^k(H)|.
		\end{align*}
		
		The last inequality follows from the fact that  $|H:A^k(T_i)|\ge7|H:A^k(H)|$ for $i=1$, $2$.
		
		Also, for $i=1$, $2$ we calculate 
		\begin{align*}
			|H:T_i||{\rm Irr}_{k,3'}(T_i)|&= \left(|H:T_i|-\frac{81}{17}\right)|{\rm Irr}_{k,3'}(T_i)|+\frac{81}{17}|{\rm Irr}_{k,3'}(T_i)|\\&\ge \frac{32}{17}|H:A^k(H)|+\frac{81}{17}|{\rm Irr}_{k,3'}(T_i)|,
		\end{align*}
		
		
		so that,
		$$\sum_{i=1}^{2} |H:T_i||{\rm Irr}_{k,3'}(T_i)|\ge \frac{64}{17}|H:A^k(H)|+\frac{81}{17}\sum_{i=1}^{2}|{\rm Irr}_{k,3'}(T_i)|.$$ 
		
		Finally, by Lemma~\ref{lem4.2}(3) we compute 
		\begin{align*}
			\displaystyle 	\sum_{\chi \in {\rm Irr}_{k,3'}(G)} \chi(1) 
			& \ge|H:A^k(H)|+5|{\rm nl}_{k,3'}(H)|+\sum_{i=1}^{l}|H:T_{i}||{\rm Irr}_{k,3'}(T_i)|\\
			& \ge|H:A^k(H)|+5|{\rm nl}_{k,3'}(H)|+\frac{64}{17}|H:A^k(H)|+\frac{81}{17}\sum_{i=1}^{l}|{\rm Irr}_{k,3'}(T_{i})|\\
			& \ge\frac{81}{17}\left(|H:A^k(H)|+|{\rm nl}_{k,3'}(H)|+\sum_{i=1}^{l}|{\rm Irr}_{k,3'}(T_{i})|\right).
		\end{align*}
		
		Applying Lemma~\ref{lem4.2}(1), we have the desired result in this case.
		
		\textbf{Case(2):} Suppose that we have at least one $H$-orbit whose length is not divisible by 3 and that length is at least 11. Among the $\alpha_{i}$, relabel the subscripts so that $|H:T_1|\ge11$. Recall that $|H:A^k(T_1)|\ge3|H:A^k(H)|$. Now, since we have that $|H:T_1|\ge11$, then we can write that $ \frac{115}{17}|H:T_1|-\frac{3(81)}{17}\ge0$. Moreover, the following inequalities hold: $3|H:T_1|-\frac{3(81)}{17}\ge\frac{64}{17}|H:T_1|$ and $|H:T_1|-\frac{81}{17}\ge\frac{64}{17(3)}|H:T_1|$.
		
		We combine these to see that 
		\begin{align*}
			\left(|H:T_1|-\frac{81}{17}\right)|{\rm Irr}_{k,3'}(T_1)|&\ge \frac{64}{3(17)}|H:T_1||T_1:A^k(T_1)|\\
			&=\frac{64}{3(17)}|H:A^k(T_1)|\\&\ge\frac{64}{17}|H:A^k(H)|.		
		\end{align*}
		
		The last inequality follows from  the fact proved in the second paragraph of this proof that $|H:A^k(T_1)|\ge3|H:A^k(H)|$. 
		
		Now, we compute 
		\begin{align*}
			|H:T_1||{\rm Irr}_{k,3'}(T_1)|& = \left(|H:T_1|-\frac{81}{17}\right)|{\rm Irr}_{k,3'}(T_1)|+\frac{81}{17}|{\rm Irr}_{k,3'}(T_i)|\\
			&\ge\frac{64}{17}|H:A^k(H)|+\frac{81}{17}|{\rm Irr}_{k,3'}(T_1)|,
		\end{align*}
		
		so that $$\sum_{i=1}^{l} |H:T_i||{\rm Irr}_{k,3'}(T_i)|\ge \frac{64}{17}|H:A^k(H)|+\frac{81}{17}\sum_{i=1}^{l}|{\rm Irr}_{k,3'}(T_i)|.$$ 
		
		Finally, by Lemma~\ref{lem4.2}(3) we compute
		\begin{align*}
			\displaystyle \sum_{\chi \in {\rm Irr}_{k,3'}(G)} \chi(1) 
			& \ge|H:A^k(H)|+5|{\rm nl}_{k,3'}(H)|+\sum_{i=1}^{l}|H:T_{i}||{\rm Irr}_{k,3'}(T_i)|\\
			& \ge|H:A^k(H)|+5|{\rm nl}_{k,3'}(H)|+\frac{64}{17}|H:A^k(H)|+\frac{81}{17}\sum_{i=1}^{l}|{\rm Irr}_{k,3'}(T_{i})|\\
			& \ge\frac{81}{17}\left(|H:A^k(H)|+|{\rm nl}_{k,3'}(H)|+\sum_{i=1}^{l}|{\rm Irr}_{k,3'}(T_{i})|\right).
		\end{align*}\\
		Applying Lemma~\ref{lem4.2}(1), we obtain the desired result.	
	\end{proof}	
\end{lemma}
\section{Main Theorem} 
In Section 5, we prove the Main Theorem. We consider a field $k$ that contains the primitive $p^{th}$ roots of unity  such that $\mathbb{Q}_p\subseteq k\subseteq\mathbb{C}$. Suppose that $t$ is the least common multiple of the orders of the roots of unity in $k$ having order dividing $|G|$. We will show that the optimal upper bound for groups of odd order depends on the odd prime $p$ and the field $k$. The argument in this proof is for groups of odd order and the reader can compare it to Theorem 4.1 of \cite{MR3552294}.

\begin{theorem} \label{thm5.1}
	Let $G$ be a group of odd order, and let $k$ be an extension of $\mathbb{Q}$ that contains the primitive $p^{th}$ roots of unity, and let $t$ be the least common multiple of the orders of roots of unity in $k$ having order dividing $|G|$. Assume one of the following:
	\begin{enumerate}
		\item $k$ contains the $3^{rd}$ roots of unity, $p\equiv1\pmod{3}$, and ${\rm acd}_k(G)<3(p+2)/(p+8)$; 
		\item $k$ contains the $3^{rd}$ roots of unity, $k$ does not contain the $({\frac{p-1}{2}})^{th}$ roots of unity or $(p-1)/2$ is even, $p\equiv2\pmod{3}$,  and ${\rm acd}_k(G)<3(p^2+2)/(p^2+8)$; 
		\item $k$ contains $({\frac{p-1}{2}})^{th}$  roots of unity, $(p-1)/2$ is odd, $p\equiv2\pmod{3}$, and ${\rm acd}_k(G)<3(p-1)/(p+3)$; 
		
		\item $k$ does not contain the $3^{rd}$ roots of unity, $p\equiv1\pmod{3}$, and ${\rm acd}_k(G)<3p/(p+2)$; 
		\item $k$ does not contain the $3^{rd}$ roots of unity, $k$ does not contain the $({\frac{p-1}{2}})^{th}$ roots of unity or $(p-1)/2$ is even, $p\equiv2\pmod{3}$, and ${\rm acd}_k(G)<3p^2/(p^2+2)$; 
		\item $p=3$,  $k$ contains the  $13^{th}$ roots of unity and ${\rm acd}_{k,3'}(G)<\frac{13}{5}$, 
		\item $p=3$, $k$  does not contain the $13^{th}$ roots of unity and ${\rm acd}_{k,3'}(G)<\frac{182}{61}$.
		
	\end{enumerate}
	Then $G$ has a normal $p$-complement. 
	\begin{proof}
		We will work by induction on $|G|$. When $G$ is an abelian group, the result is trivial. Therefore, we assume $G$ is  a nonabelian group. Thus, we have $G'>1$. Observe that $G$ is a group of odd order, so $G$ is solvable by the Feit-Thompson Theorem of \cite{MR166261}. The derived subgroup $G'$ contains a minimal normal subgroup, say $N$, of $G$, and so  all linear characters will be contained in ${\rm Irr}_{k,p'}(G/N)$. Since every character in ${\rm Irr}_{k,p'}(G)$ that is not in ${\rm Irr}_{k,p'}(G/N)$ is nonlinear, they all have degrees at least $3$, and since ${\rm acd}_{k,p'}(G)<3$, we conclude that ${\rm acd}_{k,p'}(G/N)\le{\rm acd}_{k,p'}(G)$. This implies that $G/N$ will satisfy the inductive hypothesis, and thus $G/N$ has a normal $p$-complement $K/N$. If $N$ is a $p'$-group, then $K$ will be a normal $p$-complement of $G$ and the result follows. 
		
		
		
		Suppose that $N$ is a $p$-group. Let $K_1$ be a Hall $p$-complement of $K$. Observe that $K_1$ is a Hall $p$-complement of $G$. By the Frattini argument, $G=KN_G(K_1)=NN_G(K_1)$. Since $|K:N|$ is co-prime to $p$ and $|K:K_1|$ is a power of $p$, we have $K=NK_1$. Fix $H=N_G(K_1)$.When $H=G$, $K_1$ will be the normal $p$-complement for $G$, and we are done. We assume that $H<G$. Note that $ H\cap N$ is normal in $H$ and $H\cap N$ is normal in $N$ since $N$ is abelian. Thus, $H\cap N$ is normal in $G$ and proper in $N$. This implies that $H\cap N=1$ by minimality of $N$.
		
		Now, suppose that $H$ contains a subgroup $M$ such that $M$ is a minimal normal subgroup of $G$. If $M\le G'$, then every character in ${\rm Irr}_{k,p'}(G) $ that is not in ${\rm Irr}_{k,p'}(G/M)$ has degree at least 3. Since ${\rm acd}_{k,p'}(G)<3$, we determine that ${\rm acd}_{k,p'}(G/M)<3$. Also, if $M\cap G'=1$, then by Lemma~\ref{lem2.4}, we have ${\rm acd}_{k,p'}(G/M)\le {\rm acd}_{k,p'}(G)<3$. In both cases, we see that $G/M$ satisfies the inductive hypothesis, and so, $G/M$ has a normal $p$-complement $K_1M/M$. This implies that $K_1M$ is normal in $G$. Since $K_1M \le H$, we have $N\cap K_1M=1$, and so $N$ centralizes $K_1M$.
		In particular, $N$ centralizes $K_1$, so that $N$ normalizes $K_1$. This implies that $N \subset H$, which is a contradiction. 
		Therefore, $H$ must be a core-free. 
		
		If $H$ is an abelian group, then ${\rm acd}_{k,p'}(G)={\rm acd}_{k}(G)$. In Hypothesis 1, we have that $p\equiv1\pmod3$ and $k$ contains the $3^{rd}$ roots of unity. This implies that $3$ divides $t$ and we can apply Lemma~\ref{lem3.1}(1) to get that ${\rm acd}_{k}(G)\ge3(p+2)/(p+8)$, which is a contradiction. In Hypothesis 2, we can see that $p\equiv2\pmod3$, $3$ divides $t$ and $(p-1)/2$ does not divide $t$ because $k$ contains the $3^{rd}$ roots of unity, but not the $({\frac{p-1}{2}})^{th}$ roots of unity or $(p-1)/2$ is even. Thus, we can apply Lemma~\ref{lem3.1}(2) to get that ${\rm acd}_{k}(G)\ge3(p^2+2)/(p^2+8)$, which is a contradiction. In Hypothesis 3, we can note that $p\equiv2\pmod3$, and $(p-1)/2$ divides $t$ because $k$ contains $({\frac{p-1}{2}})^{th}$ roots of unity. Apply Lemma~\ref{lem3.1}(3) to get that ${\rm acd}_k(G)\ge3(p-1)/(p+3)$, which is a contradiction. In Hypothesis 4, assume that $k$ does not contain the $3^{rd}$ roots of unity. Then, we have that $t$ is not divisible by $3$. By Lemma~\ref{lem3.1}(4) we know that  ${\rm acd}_k(G)\ge3p/(p+2)$, which is a contradiction. In Hypothesis 5, assume that $k$ does not contain the $3^{rd}$ roots of unity, and does not contain the  ${((p-1)/2)}^{th}$ roots of unity or $(p-1)/2$ is even. This implies that $t$ is not divisible by $3$ and $(p-1)/2$. Then in return, by Lemma~\ref{lem3.1}(5)  we get that ${\rm acd}_k(G)\ge3p^2/(p^2+2)$, which is a contradiction. 
		
		
		
		In Hypothesis 6, assuming that $|H|=13$, we can see that $A^k(H)=1$, because $k$ contains the $13^{th}$ roots of unity. This implies that ${\rm acd}_k(G)\ge 13/5$ by Lemma~\ref{lem3.2}, which is a contradiction. Now, assume that $|H|\ne13$. Then by Lemma~\ref{lem3.2} we get that ${\rm acd}_{k}(G)\ge182/61\ge13/5$, a contradiction. In Hypothesis 7, if $|H|=13$, then we can see that $A^k(H)\ne 1$ because $k$ does not contain the $13^{th}$ roots of unity. Thus, we can apply Lemma~\ref{lem3.2} to get that ${\rm acd}_{k}(G)\ge182/61$, a contradiction. Now, assume that $|H|\ne13$. Then by Lemma~\ref{lem3.2} we get that ${\rm acd}_{k}(G)\ge182/61$, a contradiction.  
		
		If $H$ is a nonabelian group, then Lemma~\ref{lem4.3} yields that ${\rm acd}_{k,p'}(G)\ge3$ when $p\ne3$, and ${\rm acd}_{k,p'}(G)\ge81/17$ when $p=3$ which is a contradiction to the hypothesis. Hence we are done.
	\end{proof}	
\end{theorem}


We now present the proofs of the results stated in the introduction. Conclusions~(1)-(3) of Theorem~\ref{thm1.1} follow immediately from Theorem~\ref{thm5.1}~(1)-(3) using $k=\mathbb{C}$. Conclusions~(4)-(6) in Theorem~\ref{thm1.1} follow immediately from Theorem~\ref{thm5.1}~(3)-(5) using $k=\mathbb{Q}_p$.

\medskip \medskip \textbf{Examples}: Suppose that $G$ is a Frobenius group. From Chapter 3 we know that $G$ has exactly $(|V|-1)/|H|$ nonlinear irreducible characters and these all have degree equal to $|H|$. 
\begin{enumerate}
	\item Consider $G$ of order $3p$ where $p\equiv1\pmod{3}$, $|V|=p$ and $|H|=3$. Observe that $G$ has exactly $(p-1)/3$ nonlinear irreducible characters and these all have degree equal to $3$. This implies that $G$ has exactly $3$ linear characters. Hence, ${\rm acd}(G)=3(p+2)/(p+8)$. In this case, the best bound of Theorem \ref{thm1.1}(1) is met in this group $G$. Also, with the field $k=\mathbb{Q}_{p}$ we get the best bound in Theorem \ref{thm1.1}(5) that ${\rm acd}_{\mathbb{Q}_p}(G)=3p/(p+2)$, because $G$ has exactly one linear character that has values in $k$.
	\item Assume that $G$ is of order $3p^2$ where $p\equiv2\pmod{3}$ and $(p-1)/2$ is even. In this case, we have that $|H|=3$ and $|V|=p^2$. Again, $G$ has exactly $(p^2-1)/3$ nonlinear irreducible characters and these all have degree equal to $3$. This implies that $G$ has exactly $3$ linear characters. Hence, ${\rm acd}(G)=3(p^2+2)/(p^2+8)$. Thus, the best bound of Theorem \ref{thm1.1}(2) is met with the group $G$. Furthermore, when $(p-1)/2$ is even or odd, $G$ has exactly one linear character that has values in the field $k=\mathbb{Q}_{p}$. This implies that we obtain the best bound in Theorem \ref{thm1.1}(4) and (6) that ${\rm acd}_{\mathbb{Q}_p}(G)=3p^2/(p^2+2)$.
	\item Suppose that $G$ is of order $\displaystyle \frac{1}{2}p(p-1)$ where $p\equiv2\pmod{3}$ and $(p-1)/2$ is odd. In this case, we have that $|H|= \displaystyle \frac{1}{2}(p-1)$ and $|V|=p$. Again, $G$ has exactly $2$ nonlinear irreducible characters and these all have degree equal to $\frac{1}{2}(p-1)$. This implies that $G$ has exactly $\frac{1}{2}(p-1)$ linear characters. Hence, ${\rm acd}(G)=3(p-1)/(p+3)$. In this case, the best bound of Theorem \ref{thm1.1}(3) is met with the group $G$.		
\end{enumerate}

\bibliographystyle{plain}
\bibliography{bibfile}

\begin{thebibliography}{10}

\bibitem{MR2286236}
David~S. Dummit and Richard~M. Foote.
\newblock {\em Abstract algebra}.
\newblock John Wiley \& Sons, Inc., Hoboken, NJ, third edition, 2004.

\bibitem{MR166261}
Walter Feit and John~G. Thompson.
\newblock Solvability of groups of odd order.
\newblock {\em Pacific J. Math.}, 13:775--1029, 1963.

\bibitem{MR3615444}
Nguyen~Ngoc Hung.
\newblock Characters of {$p'$}-degree and {T}hompson's character degree
  theorem.
\newblock {\em Rev. Mat. Iberoam.}, 33(1):117--138, 2017.

\bibitem{MR3096606}
I.~M. Isaacs, Maria Loukaki, and Alexander Moret\'{o}.
\newblock The average degree of an irreducible character of a finite group.
\newblock {\em Israel J. Math.}, 197(1):55--67, 2013.

\bibitem{MR1280461}
I.~Martin Isaacs.
\newblock {\em Character theory of finite groups}.
\newblock Dover Publications, Inc., New York, 1994.
\newblock Corrected reprint of the 1976 original [Academic Press, New York;
  MR0460423 (57 \#417)].

\bibitem{MR2426855}
I.~Martin Isaacs.
\newblock {\em Finite group theory}, volume~92 of {\em Graduate Studies in
  Mathematics}.
\newblock American Mathematical Society, Providence, RI, 2008.

\bibitem{MR3552294}
Mark~L. Lewis.
\newblock Variations on average character degrees and {$p$}-nilpotence.
\newblock {\em Israel J. Math.}, 215(2):749--764, 2016.

\bibitem{MR2764923}
Kay Magaard and Hung~P. Tong-Viet.
\newblock Character degree sums in finite nonsolvable groups.
\newblock {\em J. Group Theory}, 14(1):53--57, 2011.

\bibitem{MR1261638}
Olaf Manz and Thomas~R. Wolf.
\newblock {\em Representations of solvable groups}, volume 185 of {\em London
  Mathematical Society Lecture Note Series}.
\newblock Cambridge University Press, Cambridge, 1993.

\bibitem{MR3210700}
Alexander Moret\'{o} and Hung~Ngoc Nguyen.
\newblock On the average character degree of finite groups.
\newblock {\em Bull. Lond. Math. Soc.}, 46(3):454--462, 2014.

\end{thebibliography}

\end{document}